\documentclass{article}%
\usepackage{amsmath}%
\usepackage{amsfonts}%
\usepackage{amssymb}%
\usepackage{graphicx}

\begin{document}

\title{On Solving Some Trigonometric Series}
\author{Henrik Stenlund\thanks{The author is obliged to Visilab Signal Technologies for supporting this work.}
\\Visilab Signal Technologies Oy, Finland}
\date{4th August, 2013}
\maketitle

\begin{abstract}
This communication shows the track for finding a solution for a sin(kx)/k**2 series and a fresh representation for the Euler's Gamma function in terms of Riemann's Zeta function. We have found a new series expression for the logarithm as a side effect. The new series are useful both for analysis, approximations and asymptotic studies. \footnote{Visilab Report \#2013-07}
\subsection{Keywords}
infinite series, diverging series, trigonometric series, Euler's Gamma function, Riemann's zeta function
\subsection{Mathematical Classification}
Mathematics Subject Classification 2010: 11M32, 11M41, 30B50, 41A58, 33B15, 42A24, 42A32
\end{abstract}
Dedicated to my mathematics teacher Unto Salomaa.
\section{Introduction}
\subsection{General}
The following infinite series are well-known and can be found in many handbooks and tables (\cite{Wheelon1953}, \cite{Abramowitz1970}, \cite{Jolley1961}, \cite{Gradshteyn2007}). Some of these series can be solved by elementary means, some by using line integrals and residues over the complex plane \cite{Byron1970}, some requiring rather advanced methods. Laplace and Fourier transforms are useful for this purpose too.
\begin{equation}
\sum^{\infty}_{k=1}\frac{cos(k{\theta})}{k}=-ln(2sin(\frac{\theta}{2})),\     0<{\theta}<{2{\pi}} \label{eqn10}
\end{equation}
\begin{equation}
\sum^{\infty}_{k=1}\frac{cos(k{\theta})}{k^2}=\frac{\pi{^2}}{6}-\frac{\pi{\theta}}{2}+\frac{\theta^2}{4} ,\     0\leq{\theta}\leq{2{\pi}} \label{eqn12}
\end{equation}
\begin{equation}
\sum^{\infty}_{k=1}\frac{cos(k{\theta})}{k^4}=\frac{\pi{^4}}{90}-\frac{\pi^2{\theta^2}}{12}+\frac{\pi{\theta^3}}{12}-\frac{\theta^4}{48} ,\     0\leq{\theta}\leq{2{\pi}} \nonumber
\end{equation}
\begin{equation}
\sum^{\infty}_{k=1}\frac{sin(k{\theta})}{k}=\frac{(\pi-\theta)}{2},\     0<{\theta}<{2{\pi}} \nonumber
\end{equation}
\begin{equation}
\sum^{\infty}_{k=1}\frac{sin(k{\theta})}{k^3}=\frac{\pi{^2}{\theta}}{6}-\frac{\pi{\theta^2}}{4}+\frac{\theta^3}{12},\     0\leq{\theta}\leq{2{\pi}}  \nonumber
\end{equation}
\begin{equation}
\sum^{\infty}_{k=1}\frac{sin(k{\theta})}{k^5}=\frac{\pi{^4}{\theta}}{90}-\frac{\pi^2{\theta^3}}{36}+\frac{\pi{\theta^4}}{48}-\frac{\theta^5}{240},\     0\leq{\theta}\leq{2{\pi}} \nonumber
\end{equation}
There are obvious relationships between these by differentiation with respect to the parameter $\theta$. It is notable that some series are missing, like 
\begin{equation}
\sum^{\infty}_{k=1}\frac{sin(k{\theta})}{k^2} \label{eqn30}
\end{equation}
and the ones obtained by integrating it. In some instances it would be useful to transform this into a different form since its behavior is not obvious at values larger than zero. It appears that the series resists known methods and no sensible results seem to be found. In the following, we use the approach of integration, starting from equation (\ref{eqn10}). In order to succeed we need to develop some interesting expressions. Our aim is to convert the series into a form which is more suitable for various analysis and transforms.

In Section 2 we derive the intermediate relations needed to start integrating the initial equation and finish with the solution. In Section 3 we shortly present some direct consequences of the new equations, applied in various ways, like getting new series expressions for some common functions. Euler's Gamma function receives a new functional equation as well. Some simple results are exhibited which are useful while solving other infinite series. They offer closed solutions to certain parametric values. Section 4 shows what will happen if the series are truncated. Some preliminary error estimates are presented. Appendices display some graphs of selected functions, entire and truncated ones.

\section{The $sin(kx)/k^2$ Series in Terms of a Zeta Function Representation}
Before integrating equation (\ref{eqn10}), we start by differentiating the commonly known Euler's $\Gamma(x)$ in Weierstrass form, to get the Digamma function, the logarithmic derivative
\begin{equation}
\frac{1}{\Gamma(x)}=x\cdot{e^{\gamma{x}}}\prod^{\infty}_{k=1}(1+\frac{x}{k})exp({-\frac{x}{k}}), \ \ x\in{C},  \  |x|<\infty \ \nonumber
\end{equation}
\begin{equation}
ln(\Gamma(x))=-ln(x)-\gamma{x}-\sum^{\infty}_{k=1}ln[(1+\frac{x}{k})exp({-\frac{x}{k}})],  \nonumber
\end{equation}
We differentiate this with respect to $x$ and get the following
\begin{equation}
\frac{d}{dx}ln(\Gamma(x))=-\frac{1}{x}-\gamma+\sum^{\infty}_{k=1}\frac{x}{k(k+x)} \nonumber
\end{equation}
We can expand the fraction inside the sum as a binomial, swap the summations and identify it as a representation of the Riemann $\zeta(s)$ function 
\begin{equation}
\sum^{\infty}_{k=1}\frac{x}{k(k+x)}=\sum^{\infty}_{k=1}(-1)^{k+1}{x^{k}{\zeta(k+1)}} \nonumber 
\end{equation}
We can integrate back to where we started. The integration constant is zero since at $x=1$ we get the familiar Euler's gamma identity
\begin{equation}
\sum^{\infty}_{k=1}\frac{(-1)^{k+1}{\zeta(k+1)}}{k+1}=\gamma \nonumber
\end{equation}
Thus we obtain the following and mark it as $S(x)$
\begin{equation}
\sum^{\infty}_{k=1}ln[(1+\frac{x}{k})exp(\frac{-x}{k})]=\sum^{\infty}_{k=1}\frac{(-1)^k{x^{k+1}{\zeta(k+1)}}}{k+1}=S(x) \label{eqn100}
\end{equation}
We can derive this in a number of other ways too. We need another function by putting a negative argument for the same and mark that as $V(x)$
\begin{equation}
\sum^{\infty}_{k=1}ln[(1-\frac{x}{k})exp(\frac{x}{k})]=\sum^{\infty}_{k=1}\frac{(-1)^k{(-x)^{k+1}{\zeta(k+1)}}}{k+1}=V(x) \nonumber
\end{equation}
The logarithm is moved to the outside and the sum becomes
\begin{equation}
S(x)+V(x)=ln[\prod^{\infty}_{k=1}(1-\frac{x}{k})(1+\frac{x}{k})]=\sum^{\infty}_{k=1}\frac{x^{k+1}{\zeta(k+1)}((-1)^k-1)}{k+1} \nonumber
\end{equation}
Further processing gives
\begin{equation}
ln[\prod^{\infty}_{k=1}(1-\frac{x^2}{k^2})]=-\sum^{\infty}_{k=1}\frac{x^{2k}{\zeta(2k)}}{k} \nonumber
\end{equation}
This is due to disappearing of even terms in the sum and to replacing the index with a more suitable one. We recognize the infinite product as a representation of $\frac{sin(\pi{x})}{\pi{x}}$. This leads to 
\begin{equation}
ln[\frac{sin(\pi{x})}{\pi{x}}]=-\sum^{\infty}_{k=1}\frac{x^{2k}{\zeta(2k)}}{k} \label{eqn140}
\end{equation}
The range of validity for this equation is
\begin{equation}
|x|<1, \ \ x\in{C}   \nonumber
\end{equation}
We can use equation (\ref{eqn140}) while integrating equation (\ref{eqn10}), getting
\begin{equation}
-\int^{}_{}{d\theta\cdot{ln(2sin(\frac{\theta}{2}))}}=-\theta{ln(\theta)}+\theta+\sum^{\infty}_{k=1}\frac{\theta^{2k+1}{\zeta(2k)}}{(2k+1)k(2\pi)^{2k}}+C_1      \label{eqn170}
\end{equation}
and apply it in the sum as follows
\begin{equation}
\sum^{\infty}_{k=1}\frac{sin(k{\theta})}{k^2}=\theta{(1-ln(\theta))}+\theta{\sum^{\infty}_{k=1}\frac{\theta^{2k}{\zeta(2k)}}{(2k+1)k(2\pi)^{2k}}} \label{eqn180}
\end{equation}
The range of validity will be
\begin{equation}
0\leq{\theta<{2\pi}}, \ \ \theta{\in{R}}    \nonumber
\end{equation}
The integration constant $C_1$ is zero because 
\begin{equation}
\lim_{\theta\rightarrow{0}}\sum^{\infty}_{k=1}\frac{sin(k{\theta})}{k^2}={0}    \nonumber
\end{equation}
Getting $\zeta{(2k)}$ into an expression is not a sign of trouble since the values are well known, see Appendix C.
\section{Direct Consequences}
\subsection{Representations for Some Trigonometric Functions}
By differentiating equation (\ref{eqn140}) with respect to $\theta$ we will obtain
\begin{equation}
cos(\pi{x})={exp(-\sum^{\infty}_{k=1}\frac{x^{2k}{\zeta(2k)}}{k})}[1-2{\sum^{\infty}_{k=1}x^{2k}{\zeta(2k)}}] \ \label{eqn200}
\end{equation}
having a range of
\begin{equation}
|x|<1, \ \ x\in{C} \nonumber   
\end{equation}
With the equations above, we get by division
\begin{equation}
tan(\pi{x})=\frac{\pi{x}}{1-2{\sum^{\infty}_{k=1}x^{2k}{\zeta(2k)}}}  \ \label{eqn210}
\end{equation}
or inverted as
\begin{equation}
\sum^{\infty}_{k=1}\zeta(2k){x^{2k}}=\frac{1}{2}(1-\pi{x}cot(\pi{x})) \nonumber
\end{equation}
being a previously known result. 
\subsection{Connection to Euler's Gamma Function}
It is interesting to see how the Gamma function would relate to what we have found in equation (\ref{eqn140}) since it is obviously closely related to equation (\ref{eqn100}). We have the traditional functional equation 
\begin{equation}
\Gamma(1-s)\cdot{\Gamma(s)}=\frac{\pi}{sin(\pi{s})} \nonumber
\end{equation}
and insert equation (\ref{eqn140}) to it getting
\begin{equation}
\Gamma(1-s)\cdot{\Gamma(1+s)}=exp[\sum^{\infty}_{k=1}\zeta(2k)\frac{s^{2k}}{k}] \ \label{eqn310}
\end{equation}
\begin{equation}
|s|<1, \ s\in C   \nonumber
\end{equation}
This is a new functional equation for the Euler's Gamma function.
\subsection{A Series Representation for the Logarithm}
We can  read the equation (\ref{eqn310}) with the argument $s=\frac{1}{z}$ as
\begin{equation}
\Gamma(1-\frac{1}{z})\cdot{\Gamma(1+\frac{1}{z})}=exp(\sum^{\infty}_{k=1}\frac{\zeta(2k)}{k{z^{2k}}}) \ \nonumber 
\end{equation}
and with $s=1-\frac{1}{z}$ as
\begin{equation}
\Gamma(1-\frac{1}{z})\cdot{\Gamma(\frac{1}{z})}=exp(\sum^{\infty}_{k=1}\frac{{\zeta(2k)}{(1-\frac{1}{z})^{2k}}}{k}){\frac{z}{z-1}} \nonumber
\end{equation}
Taking into account the rightmost fraction, we can equate the Gammas in these two equations and after subjecting it to a logarithm, get
\begin{equation}
ln(z-1)=\sum^{\infty}_{k=1}\frac{{\zeta(2k)}{(({z-1})^{2k}-1)}}{k{z^{2k}}} \nonumber
\end{equation}
By making a yet another change of variable of
\begin{equation}
z=x+1    \nonumber
\end{equation}
we obtain
\begin{equation}
ln(x)=\sum^{\infty}_{k=1}\frac{{\zeta(2k)}{(x^{2k}-1)}}{k{(x+1)^{2k}}} \ \label{eqn440}
\end{equation}
with
\begin{equation}
x\neq{0}, \ \ x\in{C}    \nonumber
\end{equation}
This is a new series representation for the natural logarithm. The derivation of it can be made in a number of other ways as well, starting from equation (\ref{eqn310}).

\subsection{Exponential Series}
The exponential series is known as Spence function having thus far only an integral representation. We can transform it by direct substitution to the exponential function by using equations (\ref{eqn12}), (\ref{eqn180}). 
\begin{equation}
\sum^{\infty}_{k=1}\frac{exp(-k{\theta})}{k^2}=\frac{\pi{^2}}{6}-\frac{\theta^2}{4}-\theta+\theta{ln(\theta)}-\theta{\sum^{\infty}_{k=1}\frac{(-1)^k\theta^{2k}{\zeta(2k)}}{(2k+1)k(2\pi)^{2k}}} \label{eqn500}
\end{equation}
The parameter $\theta \in C$ in general and the range of validity is limited to
\begin{equation}
0\leq{Re(\theta)}, \ \  |\theta|<2\pi \nonumber   
\end{equation}

\subsection{Simple Results}
From the equations above we can get results for specific values of the parameters. Equation (\ref{eqn180}) reads at $\theta\rightarrow{2\pi}{i}$
\begin{equation}
\sum^{\infty}_{k=1}\frac{\zeta(2k)}{(2k+1)k}=ln(2\pi)-1 \nonumber
\end{equation}
From the same equation comes at $\theta=\frac{\pi}{2}$ the following 
\begin{equation}
\sum^{\infty}_{k=1}\frac{\zeta(2k)}{(2k+1)k{2^{4k}}}=ln(\frac{\pi}{2})-1+\frac{\pi}{6} \nonumber
\end{equation}
From equation (\ref{eqn140}) we can obtain at $x=\frac{1}{2}$
\begin{equation}
\sum^{\infty}_{k=1}\frac{\zeta(2k)}{k{2^{2k}}}=ln(\frac{\pi}{2})  \nonumber
\end{equation}
and at $x=\frac{1}{4}$ 
\begin{equation}
\sum^{\infty}_{k=1}\frac{\zeta(2k)}{k{4^{2k}}}=ln(\frac{\pi{\sqrt{2}}}{4})  \nonumber
\end{equation}
and at $x=\frac{3}{4}$
\begin{equation}
\sum^{\infty}_{k=1}\frac{\zeta(2k){(\frac{3}{4})^{2k}}}{k}=ln(\frac{3\pi{\sqrt{2}}}{4})  \nonumber
\end{equation}
By putting $x=\frac{1}{2}$ in equation (\ref{eqn200}) produces
\begin{equation}
\sum^{\infty}_{k=1}\frac{\zeta(2k)}{2^{2k}}=\frac{1}{2}  \nonumber
\end{equation}
With $x=\frac{1}{e}$ in equation (\ref{eqn140}) we have
\begin{equation}
\sum^{\infty}_{k=1}\frac{\zeta(2k)e^{-2k}}{k}=ln(\frac{\pi}{e})-ln(sin(\frac{\pi}{e}))  \nonumber
\end{equation}
With $x=2$ in equation (\ref{eqn440}) we reach at
\begin{equation}
\sum^{\infty}_{k=1}\frac{\zeta(2k)(2^{2k}-1)}{k3^{2k}}=ln(2)  \nonumber
\end{equation}
The results above were verified numerically too. The $\zeta(s)$ very rapidly approaches unity when $s$ grows and also some of the series are converging rapidly. We are tempted to see how accurate an estimate the first term will offer for some of these series, see Appendices A and B.

\section{Truncated Series as Approximations}
The behavior of equation (\ref{eqn500}) at small $\theta$ can be easily estimated as follows, taking only the first term of the series on the right
\begin{equation}
\sum^{\infty}_{k=1}\frac{exp(-k{\theta})}{k^2}\approx{\frac{\pi{^2}}{6}-\frac{\theta^2}{4}-\theta+\theta{ln(\theta)}-\frac{\theta^3}{84}} \nonumber
\end{equation}
It seems that the new series, equation (\ref{eqn500}), converges very fast. After seven terms in the new series ($\theta{\approx{1}}$) we get 18 decimals correct when $\theta \in R$. In the original series on the left, 30 terms are required to reach the same accuracy. Attempting to get anything like this from the equation (\ref{eqn500}) by using Taylor's series is not working due to resulting diverging series. It is interesting that on the right side, the series is not important when $\theta$ is small.

The behavior of $cos(\pi{x})$, equation (\ref{eqn200}), at small $x$ is simply estimated by taking the first term of the series
\begin{equation}
cos(\pi{x})\approx{exp[-x^2{\zeta{(2)}}]}(1-2x^2{\zeta(2)}) \nonumber
\end{equation}
It will give 12 decimals correct for the $cos(\pi{x})$ if $x<0.01$. On the other hand, if $x\approx{0.2}0$ the first term offers an accuracy of three decimals. 

The $sin(\pi{x})$ function, as opened up from equation (\ref{eqn140}),
\begin{equation}
sin(\pi{x})=\pi{x}\cdot{exp(-\sum^{\infty}_{k=1}\frac{x^{2k}{\zeta(2k)}}{k})}   \nonumber
\end{equation}
gives a surprisingly good approximation when $x<0.05$ with only the first term of the series. The accuracy is at least nine significant figures. Even at $x\approx{0.5}$ we get two significant figures correct with the first term only. Of course, the $exp()$ function can be expanded as a power series in order to simplify it. The first two terms of that would be
\begin{equation}
sin(\pi{x})\approx{\pi{x}}\cdot{(1-\sum^{\infty}_{k=1}\frac{x^{2k}{\zeta(2k)}}{k})}   \nonumber
\end{equation}
The accuracy would then be at least five significant figures with the first term of the series when $x<0.05$. At $x\approx{0.5}$ the error is about $13\%$ and does not get any better even if a lot of terms were added, due to the really crude truncation of Taylor's series of the $exp()$ function. 

\section{Discussion}
We attempted in solving equation (\ref{eqn30}) to transform it into a form which is more clearly understandable and amenable for further processing. The trigonometric series has a solution which is in a form more suitable for many analysis, in equation (\ref{eqn180}). The corresponding exponential series (\ref{eqn500}) was obtained directly from the two basic results, equations (\ref{eqn12}) and (\ref{eqn180}).

We needed to develop some intermediate equations to reach these results. Equation (\ref{eqn140}) is a fast converging series expansion for
\begin{equation}
ln(\frac{sin(\pi{x})}{\pi{x}})  \nonumber
\end{equation}
It has a $cos(\pi{x})$ counterpart in equation (\ref{eqn200}), being derived directly by differentiation. 

Since the equation (\ref{eqn180}) was derived from a starting point close to the Euler's Gamma function, it has a relationship with it, leading to equation (\ref{eqn310}). That is a nice even function of $s$.

Simple algebraic play with equation (\ref{eqn310}) brings out equation (\ref{eqn440}). It is a striking, rapidly converging series representation for the natural logarithm function. Most of the new series found are fast convergent. It suffices in many cases to take just a few terms to reach any practical accuracy. 

\newpage



\newpage

\appendix
\section{Appendix. Graphic Display of the Series for the $ln(\frac{sin(\pi{x})}{\pi{x}})$ Function}
Equation (\ref{eqn140}) is depicted in the following Figure 1 around $x=0$. It is with one, two and three terms of the series, including the exact function and error terms. 

\begin{equation}
ln[\frac{sin(\pi{x})}{\pi{x}}]=-x^{2}{\zeta(2)}-\frac{x^{4}{\zeta(4)}}{2}-\frac{x^{6}{\zeta(6)}}{3}+..., \ x\in{R} \nonumber
\end{equation}

\begin{figure}[h]
    \includegraphics[width=1.0\textwidth]{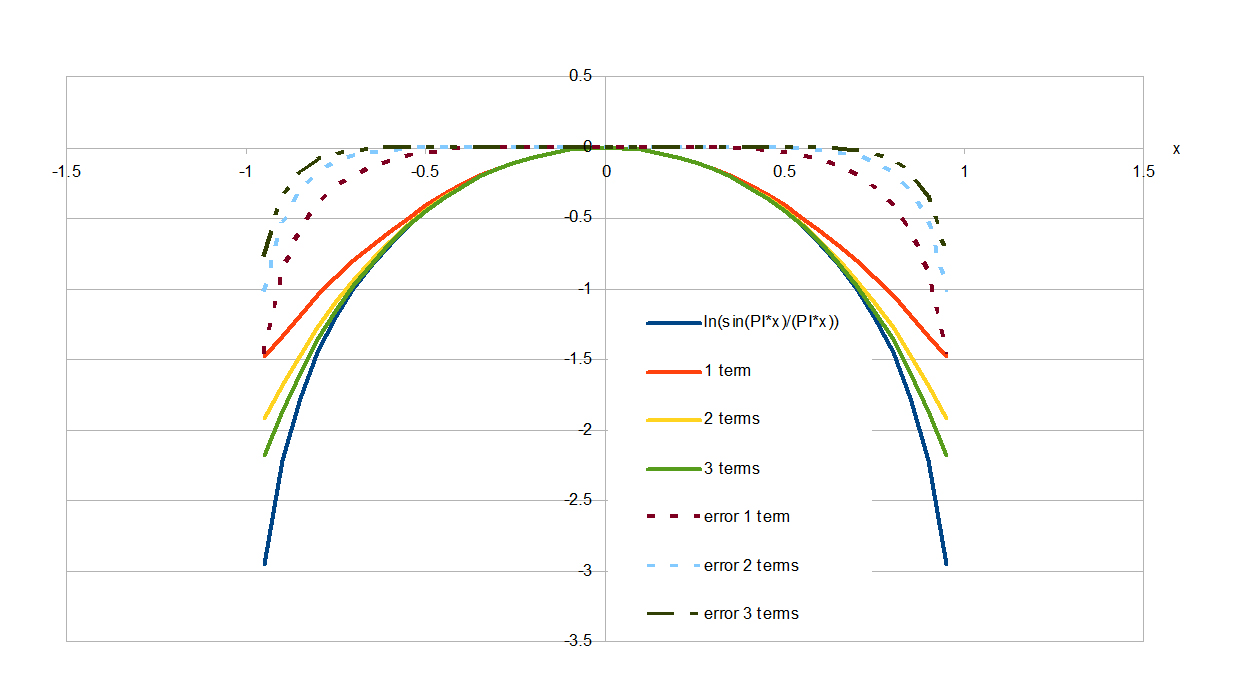}
	\label{fig:sinx}
	\caption{Logarithm of the sin(PI*x)/PI*x function approximated by one to three first terms of the series. The exact function is drawn for comparison}
\end{figure}

\newpage
\section{Appendix. Graphic Display of the Series for the Logarithm Function}
Equation (\ref{eqn440}) is drawn Figure 2 around $x=1$. It is with one, two and three terms of the series and the natural logarithm is for comparison. 

\begin{equation}
ln(x)=\frac{\zeta(2)(x^{2}-1)}{{(x+1)^{2}}}+\frac{\zeta(4)(x^{4}-1)}{{2(x+1)^{4}}}+\frac{\zeta(6)(x^{6}-1)}{{3(x+1)^{6}}}+... , \ x\in{R} \nonumber
\end{equation}

\begin{figure}[h]
    \includegraphics[width=1.0\textwidth]{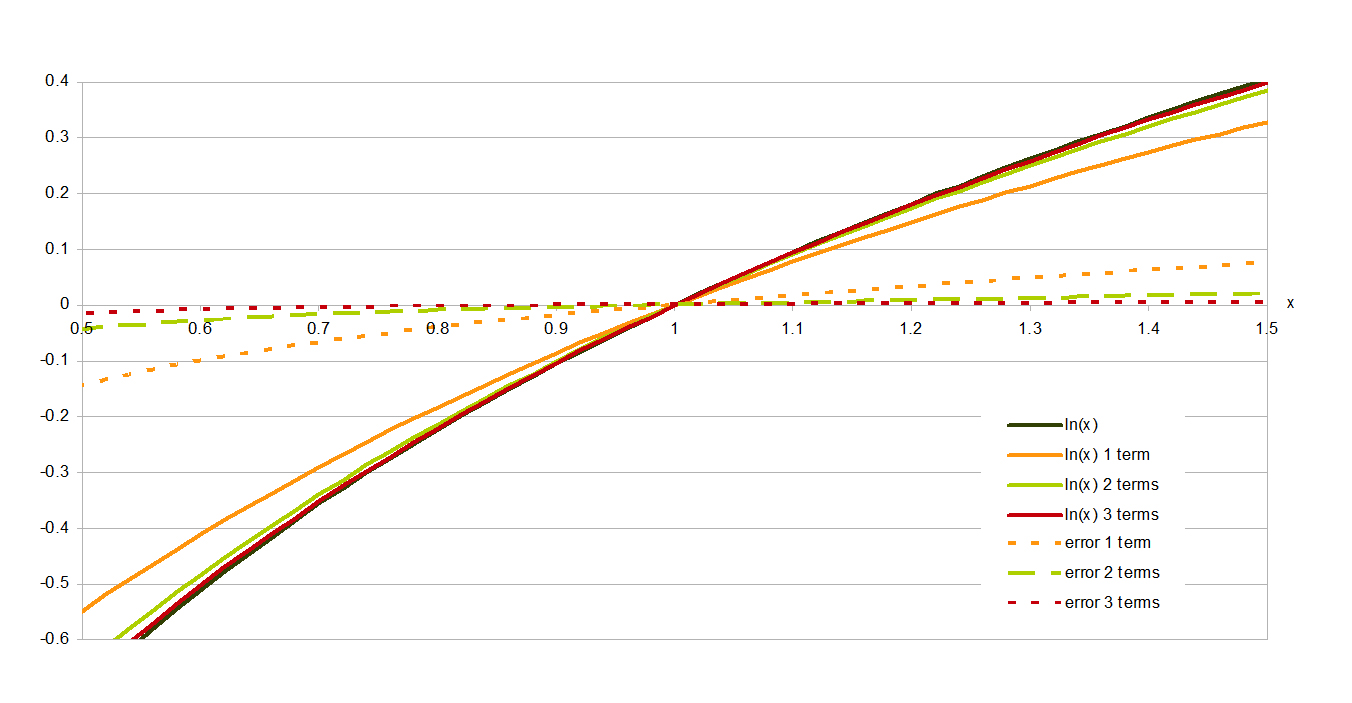}
	\label{fig:log_zeta_expansion}
	\caption{Logarithm function approximated by one to three first terms of the series. The ln(x) function is drawn for comparison}
	\end{figure}
	
\newpage
\section{Appendix. Values for the  Riemann Zeta Function for Even Integers}
In many equations above, we have the $\zeta(2k)$ function. Its values can be calculated from 
\begin{equation}
\zeta(2n)=\frac{2^{2n-1}\pi^{2n}B_n}{(2n)!}     \nonumber
\end{equation}
Here $B_n$ are the Bernoulli numbers. The first few $\zeta(2k)$ would be
\begin{equation}
\zeta(2)=\frac{\pi^{2}}{6}     \nonumber
\end{equation}
\begin{equation}
\zeta(4)=\frac{\pi^{4}}{90}     \nonumber
\end{equation}
\begin{equation}
\zeta(6)=\frac{\pi^{6}}{945}     \nonumber
\end{equation}
\begin{equation}
\zeta(8)=\frac{\pi^{8}}{9450}     \nonumber
\end{equation}
\begin{equation}
\zeta(10)=\frac{\pi^{10}}{93555}     \nonumber
\end{equation}

\end{document}